\documentclass[11pt]{article}
\usepackage[margin=4.46cm]{geometry}
\usepackage{graphicx}
\usepackage{epsfig}
\usepackage{amsmath}
\usepackage{amssymb}
\usepackage{caption}
\usepackage{tikz}
\newcommand{\commentout}[1]{}

\def \Rset {{\mathbb R}}

\def \Nset {{\mathbb N}}

\newcommand{\nit}{\noindent}

\newcommand{\be}{\begin{equation}}
\newcommand{\ee}{\end{equation}}
\newcommand{\ba}{\begin{eqnarray}}
\newcommand{\ea}{\end{eqnarray}}
\newcommand{\bi}{\begin{itemize}}
\newcommand{\ei}{\end{itemize}}
\newcommand{\br}{\begin{eqnarray}}
\newcommand{\er}{\end{eqnarray}}

\newcommand{\qed}{\mbox{$\square$}\newline}

\newtheorem{theo}{Theorem}[section]

\newtheorem{lem}{Lemma}[section]

\begin{document}
\title{Nonuniqueness of infinity ground states}
\author{ Ryan Hynd\footnote{Department of math, Courant institute, partially supported by NSF grant DMS-1004733.},\quad  Charles K. Smart\footnote{Department of math,  MIT, partially supported by NSF grant DMS-1004594.},\quad Yifeng Yu\footnote{Department of math, UC Irvine, partially supported by NSF grant DMS-0901460 and NSF CAREER award DMS-1151919.}}
\date{}
\maketitle
\begin{abstract}
In this paper,  we construct a dumbbell domain for which the associated principle $\infty$-eigenvalue is not simple.  This gives a negative answer to the outstanding problem posed in \cite{JLM}.  It remains a challenge to determine whether simplicity holds for convex domains.
\end{abstract}
\section{Introduction}
Let $\Omega$ be a bounded open set in $\Rset ^n$. According to
Juutinen-Lindqvist-Manfredi \cite{JLM}, a continuous function $u\in C(\bar \Omega)$  is said
to be an {\it infinity ground state in  $\Omega$} if it is a positive
viscosity solution of the following equation: 
\be{}
\begin{cases}
\max\left\{\lambda_{\infty}-{|Du|\over u},\ \Delta_{\infty}u
\right\}=0 & \text{in $\Omega$}\\
 u = 0 & \text{on $\partial \Omega$.}
 \end{cases}
 \ee
Here 
$$
\lambda_{\infty}=\lambda_{\infty}(\Omega)={1\over \max_{\Omega}d(x,\partial \Omega)}
$$
is the principle $\infty$-eigenvalue, and $\Delta_{\infty}$ is the infinity Laplacian operator, i.e,
$$
\Delta_{\infty}u=u_{x_i}u_{x_j}u_{x_ix_j}.
$$
The above equation is the limit as $p\to +\infty$ of the equation \be{}
\begin{cases}
-{\mathrm {div}}(|Du|^{p-2}Du)=\lambda_{p}^{p}|u|^{p-2}u &
\text{in $\Omega$}\\
u = 0 & \text{on $\partial \Omega$},
\end{cases}
\ee which is the Euler-Lagrange equation of
 the nonlinear Rayleigh quotient
$$
{\int_{\Omega}|Du|^{p}\,dx \over \int_{\Omega}|u|^{p}\,dx}.
$$
Precisely speaking, let $u_{p}$ be a positive  solution of
equation (1.2) satisfying
$$
\int_{\Omega}u^{p}_{p}\,dx=1.
$$
If $u_{\infty}$ is a limiting point of $\{u_{p}\}$, i.e, there
exists a subsequence $p_{j}\to +\infty$ such that
$$
u_{p_j}\to u_{\infty}  \quad \text{uniformly in $\bar \Omega$},
$$
it was proved in \cite{JLM} that $u_{\infty}$ is a viscosity solution
of the equation (1.1) and
$$
\lim_{p\to +\infty}\lambda_{p}=\lambda_{\infty}.
$$
We say that $u$ is {\it a variational  infinity  ground state} if it
is a limiting point of $\{u_{p}\}$.

\medskip

A natural problem regarding equation (1.1) is to deduce whether or not infinity ground states in a given domain are unique up to a multiplicative factor; in this case, $\lambda_\infty$ is said to be {\it simple}.  The simplicity of $\lambda_\infty$ has only been established for those domains where the distance function $d(x,\partial \Omega)$ is an infinity ground state (\cite{Yu}).  Such domains includes the ball, stadium, and torus.  It has been a significant outstanding open problem to verify if simplicity holds in general domains or to exhibit an example for which simplicity fails.  In this paper, we resolve this problem by constructing a planar domain where simplicity fails to hold. 

\par For $\delta\in (0,1)$, denote the dumbbell
$$
D_0=B_1(\pm 5e_1)\cup R
$$
for $R=(-5,5)\times (-\delta,\delta)$ and $e_1=(1,0)$.  Throughout this paper,  $B_r(x)$ represents the open ball centered at $x$ with radius $r$.
\begin{figure}[h]
\begin{center}
\begin{tikzpicture}
\draw[thick] (4,.2) -- (-4,.2) arc (11:349:1.02) -- (4,-.2) arc (191:360:1.02)--(6.02,0) arc (0:169:1.02);
\draw[thick,dotted] (-5,0) -- (5,0);
\fill (-5,0) circle (0.08);
\fill (5,0) circle (0.08);
\fill (0,0) circle (0.08);
\draw (0,0.5) node {$D_0$};
\draw (5,-0.3) node {$(5,0)$};
\draw (-5,-0.3) node {$(-5,0)$};
\draw (5,0) -- (5.71,0.71);
\draw (5.2,0.5) node {$1$};
\draw (-5,0) -- (-5.71,0.71);
\draw (-5.2,0.5) node {$1$};
\end{tikzpicture}
\caption{The dumbbell domain $D_0$.}
\end{center}
\end{figure}
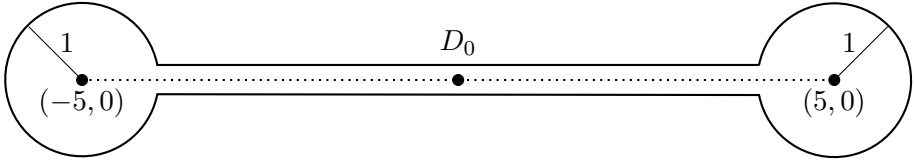
The following is our main result.
\begin{theo}\label{main} There exists $\delta_0>0$ such that when $\delta\leq \delta_0$,  the dumbbell $D_0$ possesses an infinity ground state $u_{\infty}$ which satisfies $u_{\infty}(5,0)=1$ and $u_{\infty}(-5,0)\leq {1\over 2}$.  In particular, $u$ 
is not a variational ground state and $\lambda_{\infty}(D_{0})$ is not simple.
\end{theo}

\par We remark that the infinity ground state described in the theorem is nonvariational simply because it is not symmetric with respect to the $x_2$-axis, which variational ground states can be showed to be.  This immediately follows from the fact that $\lambda_p$ is simple, which implies any solution $u_p$ of (1.2) on $\Omega=D_0$ must be symmetric with respect to the $x_2$-axis. We also remark that the number ``${1\over 2}$" in the above theorem is not special.  By choosing a suitable $\delta_0$, we can in fact make $u_{\infty}(-5,0)$ less than any positive number.  

\par For the reader's convenience, we sketch  the idea of the proof.   Consider the union of two disjoint balls with distinct radius $U_{\epsilon}=B_1(5e_1)\cap B_{1-\epsilon}(-5e_1)$ for $\epsilon\in (0,1)$.  If $u$ is an infinity ground state of $U_{\epsilon}$,  the uniqueness of $\lambda_{\infty}$ (\cite{JLM}) immediately implies that $u\equiv 0$ in $B_{1-\epsilon}(-5e_1)$.  A similar conclusion also holds for  the principle eigenfunction of $\Delta_p$. It is therefore natural to expect that such a degeneracy of $u$ on the smaller ball may change very little if we add a narrow tube connecting these two balls.  The key is to get uniform control of the width of the tube as $\epsilon\to 0$ for variational infinity ground states in an asymmetric perturbation $D_\epsilon$ of $D_0$; this is proved in Lemma \ref{keyestimate}.  Lemma \ref{keyestimate}  also implies the sensitivity of principle eigenfunctions of $\Delta_p$ when $p$ gets large.  An important step is to show that,  within the narrow tube,  the $L^p$ norm of principle eigenfunction of $\Delta_p$ is uniformly controlled by its maximum norm (Lemma \ref{control}).   We would like to point out that such a procedure as described above does not work for finite $p$.

\section{Proof}

We first prove several lemmas. Throughout this section, we write $e_1=(1,0)$ and $e_2=(0,1)$. The following estimate follows easily from comparison with the fundamental solution of the $p$-Laplacian, i.e. $|x|^{p-2\over p-1}$.
\begin{lem}\label{decay}
Let $R=(-1,1)\times (-\delta,\delta)$ for $\delta\in (0,{1\over 2})$.  Assume that $\lambda\in (0,2)$ and $u\leq 1$ is a positive solution of
\be\label{eq1}
\begin{cases}
-\Delta_pu=-\mathrm{div}(|Du|^{p-2}Du)={\lambda}^pu^{p-1} & \text{in $R$}\\
u(t,\pm \delta)=0 & \text{for $t\in [-1,1]$}.
\end{cases}
\ee
Then for $p\geq 7$
\be\label{upperBound}
u(x)\leq 6|x\pm\delta e_2|^{p-2\over p-1}.
\ee
\end{lem}
\nit Proof:  Denote  $w(x)=6|x-\delta e_2|^{\alpha}-{1\over 2}|x-\delta e_2|^{2\alpha}$ for $\alpha={p-2\over p-1}$.  Note that if $w=f(u)$, then
$$
\Delta_pw=|f'|^{p-2}f'\Delta_pu+(p-1)|f'|^{p-2}f''|Du|^p.
$$
Since $\Delta_p|x|^{\alpha}=0$ and $|x-\delta e_2|<2$, a direct computation using the above formula shows that for $p\geq 7$,
$$
-\Delta_pw=(p-1)|x-\delta e_2|^{-p\over p-1}\alpha^p(6-|x-\delta e_2|^\alpha)^{p-2}>{4^{p-3}\over 2}\geq 2^p \quad \text{in $R$}.
$$
It is straightforward to check $w>0$ in $R$ and 
$$
w(\pm 1,x_2)\geq 4  \quad \text{for $|x_2|\leq {\delta}$}.
$$
Hence
$$
u(x)\leq w(x)  \quad \text{on $\partial R$}.
$$
Combining with $-\Delta_pu\leq 2^p$,  \eqref{upperBound} follows from the comparison principle.
\qed

The following estimate may not be optimal, but is sufficient for our purposes.

\begin{lem}\label{control}
Let $R_4=(-4,4)\times (-\delta,\delta)$ for $\delta\in (0,{1\over 2})$.  Assume that $\lambda\leq 2$ and $u\leq 1$ is a positive solution of
\be\label{eq2}
\begin{cases}
-\Delta_pu=-\mathrm{div}(|Du|^{p-2}Du)=\lambda^pu^{p-1} & \text{in $R_4$}\\
u(t,\pm \delta)=0 & \text{for $t\in [-4,4]$}.
\end{cases}
\ee
Then, for $p\geq 7$ and $R_1= (-1,1)\times (-\delta,\delta)$,
\be\label{key}
\int_{R_1}u|Du|^{p-1}\,dx+\int_{R_1}|Du|^p\,dx\leq C_{0}^{p}.
\ee
Here $C_{0}>1$ is a universal constant (independent of $p$ and $\delta$).
\end{lem}

\nit Proof:  For $i=1,2,3,4$, we write $R_i=(-i,i)\times (-\delta,\delta)$.  Throughout the proof,  $C>1$ represents various numbers which are independent of $p$ and $\delta$.  We first prove an estimate which is a slight modification of a well know result (\cite{L},\cite{LM}).\par Suppose that $\xi\in C_{0}^{\infty}(R_4)$ and $0\leq \xi\leq 1$.  Multiplying  $u^{1-p}\xi^p$ on both sides of (\ref{eq2}) and using H\"older's inequality, we get 
$$
S\leq {p\over p-1}S^{1-{1\over p}}||D\xi||_{L^p(R_2)}+{2^{p+1}},
$$
where $S=\int_{R_4}|{Du\over u}|^p\xi^p\,dx$.  If ${S\over 2}\geq {2^{p+1}}$, then
$$
{S\over 2}\leq {p\over p-1}S^{1-{1\over p}}||D\xi||_{L^p(R_2)}.
$$
 Since $({p\over p-1})^p\leq 4$, we have that
\be\label{bound1}
S=\int_{R_4}\left|{Du\over u}\right|^p\xi^p\,dx\leq \max \left\{2^{p+2},\ 4\cdot 2^p\int_{R_4}|D\xi|^p\,dx \right\}.
\ee
\par Let $g_1(t)\in C_{0}^{\infty}(-4,4)$ satisfy $0\leq g_1\leq 1$, $|g_{1}^{'}|\leq 2$ and
$$
g_1(t)=1 \quad \text{for $t\in [-3,3]$}.
$$
Also, for $m\in \Nset$,  denote $\delta_m=\delta(1-{1\over m})$.  Choose $h_{m}(t)\in C_{0}^{\infty}(-\delta,\delta)$  such that $0\leq h_m\leq 1$,  $|h_{m}^{'}|\leq {2m\over \delta}$ and
$$
h_m(t)=1 \quad \text{for $t\in [-\delta_m,\delta_m$]}.
$$
For $x=(x_1,x_2)$, let $\xi_m(x_1,x_2)=g_1(x_1)h_{m}(x_2)$. Then 
$$
|D\xi_{m+1}|^p\leq 2^p(2^p+|h_{m+1}^{'}|^p)
$$
and 
$$
4\cdot 2^p\cdot \int_{R_4}|D\xi_{m+1}|^p\,dx\leq 32\cdot 8^p+32\cdot 8^{p}\cdot {\left({m+1\over \delta}\right)}^{p-1}\leq C^p\left({m\over \delta}\right)^{p-1}.
$$
Hence by (\ref{bound1})
$$
\int_{[-3,3]\times [-\delta_{m+1},\delta_{m+1}]}\left|{Du\over u}\right|^p\,dx\leq  C^p\left({m\over \delta}\right)^{p-1}.
$$
\par Owing to Lemma {\ref{decay}} and translation, we have that for $x=(x_1,x_2)\in [-3,3]\times (-\delta,\delta)$
$$
u(x_1,x_2)\leq 6\min\{(\delta-x_2)^{p-2\over p-1},\ (\delta+x_2)^{p-2\over p-1}\}.
$$
In particular, we have
$$
u(x_1,x_2)\leq 6\left({\delta\over m}\right)^{p-2\over p-1}  \quad \text{in $A_m$},
$$
where $A_m=[-3,3]\times [\delta_{m},\delta_{m+1}]$. 
Hence
$$
\int_{A_m}|Du|^p\,dx\leq C^p\cdot \left({m\over \delta}\right)^{p-1}\left({\delta\over m}\right)^{p(p-2)\over p-1}\leq C^p\cdot \left({m\over \delta}\right)^{1\over p-1};
$$
again we emphasize $C$ is independent of $p$ and $\delta$.  
\par Accordingly,
$$
\int_{[-3,3]\times [0,{\delta}]}u^2|Du|^p\,dx=\sum_{m=1}^{\infty}\int_{A_m}u^2|Du|^p\,dx\leq 36\cdot C^p\sum_{m=1}^{\infty}{1\over m^{3\over 2}}\leq C^{p}.
$$
Similarly, we can prove that
$$
\int_{[-3,3]\times [-\delta,0]}u^2|Du|^p\,dx\leq C^p,
$$
and therefore
$$
\int_{R_3}u^2|Du|^p\,dx\leq C^{p}.
$$
Using H\"older's inequality and the assumption that $u\leq 1$,  we also have that
$$
\begin{array}{ll}
\int_{R_3}u^2|Du|^{p-1}\,dx&\leq  6^{1\over p}\cdot  {(\int_{R_3}u^{2p\over p-1}|Du|^{p}\,dx)}^{p-1\over p}\\[5mm]
&\leq  2\cdot  {(\int_{R_3}u^2|Du|^{p}\,dx)}^{p-1\over p}\\[5mm]
& \leq  C^{p}.
\end{array}
$$

\par Choose $g_2(t)\in C_{0}^{\infty}(-3,3)$ such that $0\leq g_2\leq 1$, $|g_{2}^{'}|\leq 2$ and
$$
g_2(t)=1 \quad \text{for $t\in [-2,2]$}.
$$
Multiplying $w(x)=u^2\cdot g_2(x_1)$ on both sides of (\ref{eq2}) leads to
$$
\int_{R_2}u|Du|^p\,dx\leq pC^p.
$$
Again, by H\"older's inequality,  we have that
$$
\int_{R_2}u|Du|^{p-1}\,dx\leq pC^{p}.
$$
\par Finally,  select $g_3(t)\in C_{0}^{\infty}(-2,2)$ satisfying $0\leq g_3\leq 1$, $|g_{3}^{'}|\leq 2$ and
$$
g_3(t)=1 \quad \text{for $t\in [-1,1]$}.
$$
Multplying $w(x)=u\cdot g_3(x_1)$ on both sides of (\ref{eq2}) leads to
$$
\int_{R_1}|Du|^p\,dx\leq p^2C^p.
$$
Since $3^p>p^2$,  we have that
$$
\int_{R_2}u|Du|^{p-1}\,dx+\int_{R_1}|Du|^p\,dx\leq 2p^2C^p\leq (6C)^p=C_{0}^{p}.
$$
Consequently, (\ref{key}) holds, as desired.
\qed

\medskip

\begin{figure}[h]
\begin{center}
\begin{tikzpicture}
\draw[thick] (4,.2) -- (-4.2,.2) arc (14:346:.83) -- (4,-.2) arc (191:360:1.02)--(6.02,0) arc (0:169:1.02);
\draw[thick,dotted] (-5,0) -- (5,0);
\fill (-5,0) circle (0.08);
\fill (5,0) circle (0.08);
\fill (0,0) circle (0.08);
\draw (0,0.5) node {$D_\epsilon$};
\draw (5,-0.3) node {$(5,0)$};
\draw (-5,-0.3) node {$(-5,0)$};
\draw (5,0) -- (5.71,0.71);
\draw (5.2,0.5) node {$1$};
\draw (-5,0) -- (-5.58,0.58);
\draw (-5,0.5) node {$1-\epsilon$};
\end{tikzpicture}
\caption{The asymetric dumbbell domain $D_{\epsilon}$}
\end{center}
\end{figure}
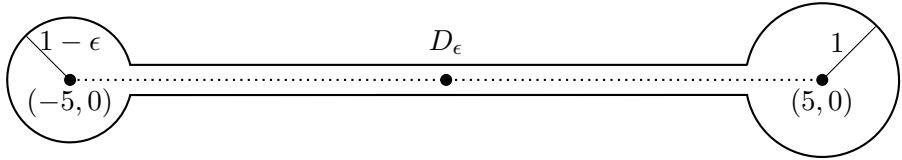 

\par Now let
$$
\delta_0={1\over 16C_0}<{1\over 16}.
$$
Here $C_0>1$ is the same number in Lemma \ref{control}.  For $\epsilon\in (0,{1\over 2})$, write
$$
D_{\epsilon}=B_{1-\epsilon}(-5e_1)\cup R\cup B_{1}(5e_1)
$$
and $R=(-5,5)\times (-\delta,\delta)$. Note that $D_{\epsilon}$ is not symmetric with respect to the $x_2$ axis and $\max_{D_{\epsilon}}d(x,\partial D_{\epsilon})=1$.

\medskip

The following lemma says that the principle eigenfunction of p-Laplacian, although unique up to multiplicative factor, is actually very sensitive to the domain when $p$ gets large.

\begin{lem}\label{keyestimate} Assume $0<\epsilon<{1\over 2}$. If $\delta\leq \delta_0$ and $u_{\infty}$ is a variational infinity ground state of $D_{\epsilon}$ satisfying $u_{\infty}(5,0)=1$, then
$$
u_{\infty}(-5,0)<{1\over 2}.
$$
Note that $\delta_0$ is independent of $\epsilon$. 
\end{lem}
\nit Proof: We argue by contradiction and assume that $u_{\infty}(-5,0)\geq {1\over 2}$.  Now fix $\delta$ and $\epsilon$. Since $\max_{D_{\epsilon}} u_{\infty}=u_{\infty}(5,0)=1$, $u_{\infty}(x)\leq d(x,\partial D_{\epsilon})$(\cite{JLM}). Hence
$$
u_{\infty}\leq \delta\leq  \delta_0   \quad \text{in $[-4,4]\times [-\delta,\delta]$}.
$$
For $p>2$, let $u_p$ be the principle eigenfunction of $\Delta_p$ in $D_\epsilon$ satisfying $\max_{D_{\epsilon}}u_p=1$ and
\be\label{eq3}
-\Delta_pu_p=-\mathrm{div}(|Du_p|^{p-2}Du_p)=\lambda_{\epsilon,p}^{p}u_{p}^{p-1} \quad \text{in $D_{\epsilon}$}.
\ee
Here $\lambda_{\epsilon,p}$ is the principle eigenvalue of $\Delta_p$ associated with $D_{\epsilon}$. 

\par Passing to a subsequence if necesary,  we may assume that
$$
\lim_{p\to +\infty}u_p=u_{\infty}  \quad \text{uniformly in $D_{\epsilon}$}.
$$
Hence, when $p$ is large enough,
\be\label{small}
u_p\leq 2\delta_0   \quad \text{in $[-4,4]\times [-\delta,\delta]$}.
\ee
Since $\lim_{p\to +\infty}\lambda_{\epsilon,p}=\lambda_{\epsilon,\infty}=1$, we may assume that $\lambda_{\epsilon,p}\leq 2$.

\par Now, define $g(t)$ by
$$
\begin{cases}
g(t)=1 & \text{for $t\leq -1$}\\
g(t)={1\over 2}(1-t) & \text{for $-1\leq t\leq 1$}\\
g(t)=0 & \text{for $t\geq 1$}.
\end{cases}
$$
Let
$$
w(x)=u_p\cdot g(x_1).
$$
and, for $\tilde R=(-5,4)\times (-\delta,\delta)$, let
$$
\Omega_{\epsilon}=B_{1-\epsilon}(-5e_1)\cup \tilde R.
$$

\begin{center}
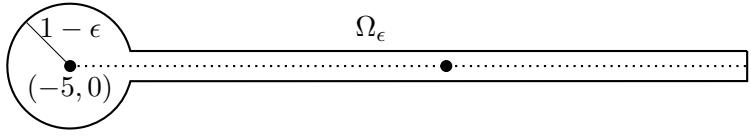

\begin{tikzpicture}
\draw[thick] (4,.2) -- (-4.2,.2) arc (14:346:.83) -- (4,-.2) -- (4,.2);
\draw[thick,dotted] (-5,0) -- (4,0);
\fill (-5,0) circle (0.08);
\fill (0,0) circle (0.08);
\draw (-1,0.5) node {$\Omega_\epsilon$};
\draw (-5,-0.3) node {$(-5,0)$};
\draw (-5,0) -- (-5.58,0.58);
\draw (-5,0.5) node {$1-\epsilon$};
\end{tikzpicture}
\captionof{figure}{The domain $\Omega_{\epsilon}$}
\end{center}
Note that  $\{w\ne 0\}\subset \Omega_{\epsilon}$ and therefore
\be\label{oneside}
\Lambda_{\epsilon,p}^{p}\leq {\int_{\Omega_{\epsilon}}|Dw|^{p}\,dx\over \int_{\Omega_{\epsilon}}|w|^{p}\,dx}={\int_{D_{\epsilon}}|Dw|^{p}\,dx\over \int_{D_{\epsilon}}|w|^{p}\,dx},
\ee
where $\Lambda_{\epsilon,p}$ is the principle eigenvalue of $\Delta_p$ associated with $\Omega_{\epsilon}$.

Since $u_p$ is uniformly H\"older continuous and $\lim_{p\to +\infty}u_p(-5e_1)=u_{\infty}(-5e_1)$, there exists $\tau\in (0,1)$ such that
\be\label{large}
u_p(x)\geq {1\over 3}  \quad \text{in $B_{\tau}(-5e_1)$},
\ee
for sufficiently large $p$.

To simplify notation, we now drop the $p$ dependence and write $u_p=u$.   Multiplying $ug^p(x_1)$ on both sides of (\ref{eq3}), we have that
$$
{\int_{D_{\epsilon}}|Du|^{p}g^{p}\,dx\over \int_{D_{\epsilon}}|w|^{p}\,dx}\leq \lambda_{\epsilon,p}^{p}+{p\int_{[-1,1]\times [-\delta,\delta]}u|Du|^{p-1}\,dx\over \int_{D_{\epsilon}}|w|^{p}}.
$$
Due to Lemma \ref{control} and (\ref{small})
$$
\int_{[-1,1]\times [-\delta,\delta]}u|Du|^{p-1}\,dx\leq (2\delta_0C_0)^p< {1\over 4^p}.
$$
Therefore owing to (\ref{large}),
$$
{p\int_{[-1,1]\times [-\delta,\delta]}u|Du|^{p-1}\,dx\over \int_{D_{\epsilon}}|w|^{p}}\leq \left({3\over 4}\right)^{p}{p\over \pi \tau^2}.
$$
\par Since $Dw=gDu+uDg$ and $(a+b)^p\leq 2^p(a^p+b^p)$, we have that
$$
\begin{array}{ll}
\int_{D_{\epsilon}}|Dw|^p\,dx&\leq \int_{D_{\epsilon}}|Du|^pg^p\,dx+2^p\int_{[-1,1]\times [-\delta,\delta]}(|Du|^pg^p+{u^p\over 2^p})\,dx\\[6mm]
&\leq \int_{D_{\epsilon}}|Du|^pg^p\,dx+(\delta_04C_0)^p+(2\delta_0)^p\\[6mm]
&\leq \int_{D_{\epsilon}}|Du|^pg^p\,dx+ 2\cdot {1\over 4^p}.
\end{array}
$$
The first inequality is also due to the fact that 
$$
Dw=gDu    \quad \text{in $D_{\epsilon}\backslash [-1,1]\times [-\delta,\delta]$ }.
$$
Therefore by (\ref{large})  when $p$ is large enough
\be\label{otherside}
{\int_{D_{\epsilon}}|Dw|^{p}\,dx\over \int_{D_{\epsilon}}|w|^{p}\,dx}\leq\lambda_{\epsilon,p}^{p}+3\cdot \left({3\over 4}\right)^{p}{p\over \pi \tau^2}\leq \lambda_{\epsilon,p}^{p}+1.
\ee
Since $\max_{D_\epsilon} d(x,\partial D_\epsilon) = 1$ and $\max_{\Omega_\epsilon} d(x,\partial \Omega_\epsilon) = 1 - \epsilon$, we have $\Lambda_{\epsilon,p} \to (1-\epsilon)^{-1}$ and $\lambda_{\epsilon, p} \to 1$ as $p \to \infty$. 

\par Thus, for sufficiently large $p$, we have
$$
\Lambda_{\epsilon,p}\geq {1\over 1-{1\over 2}\epsilon}\quad \mathrm{and}\quad \lambda_{\epsilon,p}\leq {1\over 1-{1\over 4}\epsilon}.
$$
Owing to (\ref{oneside}) and (\ref{otherside}), we have
$$
\left({2\over 2-\epsilon}\right)^{p}\leq \left({4\over 4-\epsilon}\right)^p+1,
$$
for all large enough $p$. Since this is a contradiction, the lemma follows. \qed

\medskip

\nit{\bf Proof of Theorem \ref{main}}:  For $\epsilon\in (0,{1\over 2})$, let $u_{\epsilon,\infty}$ be a variational infinity ground state of $D_{\epsilon}$ satisfying $u_{\epsilon,\infty}(5,0)=1$.  Since $\Delta_{\infty}u_{\epsilon,\infty}\leq 0$, according to \cite{CEG},  the sequence $\{u_{\epsilon,\infty}\}_{\epsilon>0}$ is uniformly Lipschitz continuous within any compact subset of $D_0$ when $\epsilon$ is small.  It is also controlled by $0\leq u_{\epsilon,\infty}\leq d(x, \partial D_{\epsilon})$ near the boundary.  Upon a subsequence if necessary, we may assume that
$$
\lim_{\epsilon\to 0}u_{\epsilon,\infty}=u_{\infty}.
$$
Then according to Lemma \ref{keyestimate},  $u_{\infty}$ is an infinity ground state of $D_0$ satisfying
$$
u_{\infty}(-5,0)\leq {1\over 2}\quad \mathrm{and}\quad  u_{\infty}(5,0)=1.
$$
As $u_{\infty}$ is not symmetric about the $x_2$-axis, it cannot be a variational infinity ground state associated to $D_0$.  As there exists at least one variational ground state \cite{JLM}, it follows that $\lambda_{\infty}(D_0)$ is not simple.

\bibliographystyle{plain}

\end{document}